\documentclass{amsart}

\usepackage[all] {xy}
\usepackage[dvips]{graphicx}
\usepackage{graphics,color}
\usepackage{tikz}
\usepackage{amssymb}
\title[Symbolic Extensions and Dominated Splittings]{Symbolic Extensions and dominated splittings for Generic $C^1$-Diffeomorphisms}

\author[A. Arbieto]{A. Arbieto}
\address{Instituto de Matem\'atica, Universidade Federal do Rio de Janeiro, P. O. Box 68530, 21945-970 Rio de Janeiro, Brazil.}
\email{arbieto@im.ufrj.br}

\author[A. Armijo ]{A. Armijo}
\address{Instituto de Matem\'atica, Universidade Federal do Rio de Janeiro, P. O. Box 68530, 21945-970 Rio de Janeiro, Brazil.}
\email{almaar@im.ufrj.br}

\author[T. Catalan]{T. Catalan}
\address{Faculdade de Matem\'atica, Universidade Federal de Uberl\^andia.}
\email{tcatalan@famat.ufu.br}

\author[L. Senos ]{L. Senos}
\address{Instituto de Matem\'atica, Universidade Federal do Rio de Janeiro, P. O. Box 68530, 21945-970 Rio de Janeiro, Brazil.}
\email{senos@im.ufrj.br}

\thanks{Partially supported by CAPES, CNPq, FAPERJ (``Jovem Cientista do Nosso Estado'') and PRONEX/DS from Brazil.}

\newtheorem{theorem}{Theorem}
\newtheorem{corollary}[theorem]{Corollary}

\newtheorem{lemma}[theorem]{Lemma}
\newtheorem{proposition}[theorem]{Proposition}

\theoremstyle{definition}
\newtheorem{definition}[theorem]{Definition}
\newtheorem{remark}[theorem]{Remark}

      \newcommand{\La}{\Lambda}

\newcommand{\mundo}{\operatorname{Diff}^1(M)}

\newcommand{\Z}{\mathbb{Z}}
\newcommand{\N}{\mathbb{N}}
\newcommand{\R}{\mathbb{R}}

\newcommand{\eps}{\varepsilon}

\newcommand{\SM}{{\mathcal M}}

\newcommand{\SR}{{\mathcal R}}
\renewcommand{\SS}{{\mathcal S}}

\newcommand{\SU}{{\mathcal U}}

\newcommand{\diff}{\operatorname{Diff}}

\begin{document}

\begin{abstract}
Let $\mundo$ be the set of all $C^1$-diffeomorphisms $f:M\rightarrow M$, where $M$ is a compact boundaryless d-dimensional manifold, $d\geq2$.
We prove that  there is a residual subset $\mathfrak R$ of $\mundo$ such that if $f\in\mathfrak R$ and if $H(p)$ is the homoclinic class associated
with a hyperbolic periodic point $p$, then either $H(p)$ admits a dominated splitting of the form $E\oplus F_1\oplus \dots\oplus F_k\oplus G$, where $F_i$ is not hyperbolic and one-dimensional, or $f|_{H(p)}$ has no symbolic extensions.
\end{abstract}

\maketitle

\section{Introduction}

\emph{Expansiveness} is an important notion in the theory of dynamical systems. Let $M$ be a compact manifold and $f:M\to M$ be a homeomorphism. Roughly, it says that if the orbits of different points must separate in finite time.  More precisely, there exists $\eps>0$ such that for any point $x$, the $\eps$-set of $x$, given by the points $y$ such that $d(f^n(x),f^n(y))<\eps$ for every integer $n$, reduces to the point $x$. This notion is somewhat related with the well know notion of sensitivity to initial conditions, commonly known as \emph{chaos}, which means that for any point, there exists a point such that the future orbit of these two points separated.  Moreover, expansiveness naturally appears in hyperbolic sets, and together with the shadowing property, play a central role to prove their stability.

However, it is important to look for weaker forms of expansivity. Clearly, expansivity implies \emph{h-expansivity}, i.e. for some $\eps>0$  the entropy of the $\eps$-set of any point $x$ is zero. This notion implies semicontinuity of the entropy map, henceforth leading to the existence of equilibrium states, which is a well know problem in ergodic theory. We remark that $h$-expansiveness do not imply expansivity. This weaker property holds for partially hyperbolic diffeomorphisms such that their central subbundle admits a dominated splitting by one-dimensional subbundles, see \cite{D-F-P-V}. It also holds for diffeomorphisms away from tangencies see \cite{L-V-Y}. 

It turns out that $h$-expansiveness implies the existence of \emph{symbolic extensions}, see \cite{B-F-F}. This means that the system is a quotient of a subshift of finite type. Actually, we can ask if the residual entropy of this extension is zero, in this case we say that the extension is \emph{principal}. However, the existence of symbolic extensions does not imply any kind of expansiveness, even asymptotic $h$-expansiveness, which requires that the entropy of the $\eps$-sets goes to zero if $\eps$ goes to zero.  In particular, \emph{the non existence} of symbolic extensions implies that a positive amount of entropy, far from zero, can be found in arbitrarily small sets, given some complexity of the dynamics, see \cite{BD}.  In the other hand, symbolic extensions are used in the theory of data transmission, see \cite{Down}. It is worthing to remark that any $C^{\infty}$ diffeomorphism is asymptotically $h$-expansive, see \cite{Bu}.

The existence of symbolic extensions is somewhat rare in non-hyperbolic dynamics. Indeed, it was proved by \cite{D-N} that $C^1$-generic non-Anosov symplectic diffeomorphisms in surfaces do not have symbolic extensions. This result was extended to higher dimensions by \cite{C-T}. By generic, we mean that this holds for systems in a residual subset of such diffeomorphisms.

A natural question in dynamical systems is to know whether the presence of a dynamical property in a $C^1$-robust way implies some hyperbolicity. For instance, \cite{B-D-P} shows that robust transitivity implies the existence of a dominated splitting. Naturally, some authors asked this question using expansivity. Indeed, Ma\~n\'e \cite{Mane} shows that any robustly expansive diffeomorphisms is Axiom A. The same question can be asked in a semi-local way. More precisely, we can ask if a homoclinic class has some expansiveness in a robust way then it is hyperbolic. By homoclinic class we mean the closure of the transversal homoclinic intersections of a periodic orbit. The series of papers \cite{PPV}, \cite{PPSV}, \cite {SV1} and \cite{SV2}, essentially proves that robustly expansive homoclinic classes are hyperbolic, see the articles for more details.  In \cite{P-V}, it was proved that any robustly $h$-expansive homoclinic class has a dominated splitting of the form $E\oplus F_1\oplus \dots\oplus F_k\oplus G$, where $F_i$ is not hyperbolic and one-dimensional. A related result was proved by Li in the context of R-robustly h-expansive homoclinic classes, see \cite{LI} for more details. %[ABCDW] !!!!! AQUI ACHO Q NÃO PRECISA ENTRAR ABCDW, AFINAL ESTAMOS FALANDO DA IMPLICAÇÃO DE PROPRIEDADES ROBUSTAS! DEIXA O BICHO PARA ENTRAR LÁ EM BAIXO! :) 

Another related question is the existence of a residual subset where the presence of a dynamical property implies hyperbolicity in the global and semi-local case. For instance, in \cite{Ar} it is proved that any generic expansive diffeomorphism is Axiom A. In \cite{C}, it was proved that generic volume preserving diffeomorphisms have symbolic extensions if, and only if, they are partially hyperbolic. In the semi-local case,  \cite{GY} proved that for a generic diffeomorphisms, any expansive homoclinic class is hyperbolic.

In this article we study these questions for generic diffeomorphisms in the semi-local case but using symbolic extensions, that as we saw before, is much weaker than expansiveness. Another results dealing with the non-existence of symbolic extensions are: \cite{D-F} constructed a locally residual subset of $C^1$-partially hyperbolic diffeomorphisms without symbolic extensions, \cite{A} also constructed other examples, for smoother systems \cite{D-N} conjectured that $C^r$-diffeomorphisms have symbolic extensions if $r\geq 2$, \cite{B} proved this conjecture for surfaces diffeomorphisms, \cite{BF} extended this result for higher dimensions with 2-dimensional center subbundle. Any $C^r$-one-dimensional transformation, with $r>1$, has symbolic extensions, this was proved by \cite{DM}.

Now, we give precise definitions and state our main results.

We consider a compact boundaryless $d$-dimensional Riemmanian manifold $M$, $d\geq2$, and denoting by $\diff^r(M)$ the set of $C^r$ diffeomorphisms on $M$ endowed with the $C^r$ topology.

\begin{definition}
A dynamical system $f:M\to M$ has a \textit{symbolic extension} if there exists a subshift $\sigma:N\to N$ and a continuous surjective map
$\pi:N\rightarrow M$ such that $\pi\circ \sigma = f \circ \pi$. In this case the system $\sigma:N\to N$ is called an \textit{extension} of $f:M\to M$ and $f$
 is called a \textit{factor} of $\sigma$. If $h_{\pi^*\mu}(f)=h_{\mu}(\sigma)$ for every invariant measure $\mu$ of $\sigma$ then the extension is called \emph{principal}.
\end{definition}

We say that $f:M^n\to M^n$ has a good decomposition in $\Lambda$ if there exists a dominated splitting $T_{\Lambda}M=E_1\oplus\dots\oplus E_k$ such that $dim(E_1)=s$, $dim(E_k)=n-u$ and for every $1<j<k$ we have $dim(E_j)=1$. Here, $s$ (resp. $u$) denotes the smallest (resp. greatest) index of a hyperbolic periodic point in $\Lambda$. Recall that  the {\it index} of a hyperbolic periodic point  $p$ is the dimension of its stable manifold.

\begin{theorem}\label{maintheo}
There is a residual subset $\mathcal{R}$ of $\mundo$ such that if $f\in\mathcal{R}$, then for every homoclinic class $H(p,f)$,
\begin{itemize}
 \item[a)] either $H(p,f)$ has a good decomposition,
\item[b)]  or $f|_{H(p,f)}$  has no symbolic extensions.
\end{itemize}
\end{theorem}

To prove this theorem we will use criterions to the non existence of symbolic extensions developed by  Downarowicz and Newhouse in \cite{D-N}.  Moreover, we also use a dichotomy between good decompositions and the existence of a homoclinic tangency, see \cite{A-B-C-D-W}. %LEMBRAR QUE FALOU ANTES.

Even so, once that one obtain a good decomposition, is somewhat folklore to obtain partial hyperbolicity when the class is isolated. In particular, we obtain the following theorem and prove it just for sake of completeness.

\begin{theorem}
There is a residual subset $\mathcal{R}$ of $\mundo$ such that if $f\in\mathcal{R}$, then for every isolated homoclinic class $H(p,f)$
\begin{itemize}
 \item[a)] either $H(p,f)$ is partially hyperbolic,
\item[b)]  or $f|_{H(p,f)}$ has no symbolic extensions.
\end{itemize}
\end{theorem}

However, this theorem together with the result of Diaz, Fisher, Pac\'ifico and Vieitez \cite{D-F-P-V}  has an interesting directly consequence.

\begin{corollary}
There is a residual subset $\mathcal{R}$ of $\mundo$ such that if $f\in\mathcal{R}$, any isolated homoclinic class of $f$ has a symbolic extension if, and only if, it has a principal symbolic extension.
\end{corollary}

Finally, as a byproduct of the techniques used in the proof of the main theorem we also get the following interesting consequence, which is somewhat related to the previous result by Pac\'ifico and Vieitez   \cite{P-V} mentioned before.

\begin{proposition}
Let $HT\subset\mundo$ be the set of diffeomorphisms exhibiting a homoclinic tangency, and  $NAHE\subset\mundo$ be  the set of  diffeomorphisms that are not asymptotically $h-$expansive.  Then $\overline{HT}=\overline{NAHE}$.
\label{prophexp}\end{proposition}

As a consequence, if a diffeomorphism is stably asymptotically $h$-expansive then it has a dominated splitting in the pre-periodic set, using a result of Wen, see \cite{W}. Moreover, if the diffeomorphism is generic then it is partially hyperbolic due to \cite{C-S-Y}.

This article is organized as follows: In Section \ref{definition}, we define precisely the notions and objects used in this paper, in Section \ref{SNP} we define and study the $\SS_{n,p}$ property, which is our tool to find diffeomorphisms that has no symbolic extensions, in Section \ref{localversion} we prove a local version of the Theorem \ref{maintheo}, in Section \ref{mainproof} we give a proof for Theorem \ref{maintheo}, in Section \ref{isolated} we consider the isolated case and, finally, in Section \ref{hexpansivity} we prove Proposition \ref{prophexp}.

\section{Definitions}\label{definition}

In this section we define precisely the notions and objects used in the introduction.

We say that $p$ is a \emph{periodic point} if $f^n(p)=p$ for some $n\geq 1$, the minimal such natural is called the \emph{period} of $p$ and it is denoted by $\tau(p,f)$, or simply by $\tau(p)$ if the diffeomorphisms $f$ is fixed. The periodic point is \emph{hyperbolic} if the eigenvalues of $Df^{\tau(p)}(p)$ do not belong to $S^1$.

If $p$ is a hyperbolic periodic point then its \emph{homoclinic class} $H(p,f)$ is the closure of the transversal intersections of the stable manifold and unstable manifold of the orbit of $p$:
$$H(p,f)=\overline{W^s(p)\pitchfork W^u(p)}.$$
It is well known that a homoclinic class is transitive. Moreover, we say that a hyperbolic periodic point $q$ is \emph{related} to $p$ if $W^s(p)\pitchfork W^u(q)\neq \emptyset$ and $W^u(p)\pitchfork W^s(q)\neq \emptyset$, it can be proved that the homoclinic class of $p$ is also the closure of the hyperbolic periodic points related to $p$.

Let $Per^n_h(f)$ be the collection of hyperbolic periodic points of $f$ of period less than or equal to $n$, and let $Per_h(f)=
{\displaystyle\bigcup_{n\geq 1} Per^n_h(f)}$.

\subsection{Domination}

%%%%%%%%%%%%%%%%%%%%%%%%%%%%%%%%%%%%%%%%%%%%%%%%%%%
We say that a compact $f$-invariant set $\La \subset M$ admits a \emph{dominated splitting}
if the tangent bundle $T_\La M$ has a continuous $Df$-invariant splitting $E_1\oplus \cdots \oplus E_k$ and there exist constants  $C > 0$,
 $0 < \lambda < 1$, such that
\begin{equation*}
||Df^n|{E_i(x)}||\cdot ||Df^{-n}|E_j(f^n(x))||\leq C\lambda^n,\quad\forall x\in\Lambda,\; n\geq0, \text{ for every } i<j.
\end{equation*}

We say that  $T_{\Lambda}M=E_1\oplus\ldots\oplus E_k$ is the finest dominated splitting if there is no dominated splitting of $E_l$ for every $1<l<k$.

%\end{definition}
%In other
%words, $\Lambda$ has a good decomposition if the central bundle admits a decompostion in unidimensional sub-bundles.

%\begin{proposition}
% Suppose $\xi$ is a compact subset of $\SM( f )$ such that there is a positive real number $\rho_0$ such that for each %$\mu \in \xi$ and each
%$k > 0$,
%$$\limsup_{\nu \in \xi,\nu \rightarrow \mu} [h_\nu(f) - h_k(\nu)] > \rho_0.$$
%Then, $f$ admits no symbolic extensions.
%\label{propext}\end{proposition}

\subsection{Hyperbolicity}

If $\La$ is a compact invariant set of a diffeomorphism $f$ then $\La$ is said to be a \textit{hyperbolic set} if we have a $Df$-invariant continuous splitting $T_{\Lambda}M=E^s\oplus E^u$ and constants $C > 0$ and $\kappa < 1$ such that
$$||Df^{-n}(x)_{|E^u_x}||\leq C\kappa^n\mbox{ and }||Df^n(x)_{|E^s_x}||\leq  C\kappa^n,$$ for every $x \in \La$ and $n \in \N$.

Let $E\oplus F_1\oplus\dots\oplus G$ be a dominated splitting over $\Lambda$. If $E$ contracts and $G$ expands, like in the previous paragraph then we say that $\Lambda$ is \emph{partially hyperbolic}.

Let $\La$ be a hyperbolic set for $f$. We call $\La$ a \textit{hyperbolic basic set} if
\begin{itemize}
 \item it is \emph{isolated}, i.e. there is a neighborhood $U$ of $\La$ such that $$\bigcap_{n \in \Z} f^n(U) = \La\,\ \textsl{and}$$
\item  $f$ has a dense orbit in $\La$.
\end{itemize}

\subsection{Genericity}

We say  that a subset $\mathfrak R\subset \mundo$ is a \textit{residual subset} if contains a countable intersection of open and dense sets.

The countable intersection of residual subsets is also a residual subset.
Since $\mundo$ is a Baire space when endowed with the $C^1$-topology, any residual subset of $\mundo$ is dense.

We will say that a property (P) holds \emph{generically} if there exists a residual subset $\mathfrak R$ such that any $f\in \mathfrak R$ has the
property (P).

\subsection{Measures and Exponents}

A measure $\mu$ is $f$-invariant if $\mu(f^{-1}(B))=\mu(B)$ for every measurable set $B$. An invariant measure is ergodic if the measure of any invariant set is zero or one. Let $\SM(f)$ be the space of $f$-invariant \textit{probability measures} on $M$, and let $\SM_e(f)$ denote the ergodic elements of $\SM(f)$.

For a hyperbolic periodic point $p$ of $f$ with period $\tau(p)$, we let $ \mu_p$ denote the periodic measure given by
$$\mu_p =\frac{1}{\tau(p)} \sum_{ x\in O(p)} \delta_x $$
where $O(p)$ denotes the orbit of $p$ and $\delta_x$ is the Dirac measure at $x$.

A measure $\mu \in \SM(f)$ is called a \textit{hyperbolic measure} for $f$ if its topological support $supp(\mu)$ is contained in a hyperbolic
basic set for $f$.
\footnote{This differs from usual definition, where one ask that the Lyapunov exponents of the measure are non-zero.}

Let $C(M,\R)$ be the set of all continuous functions $h:M\rightarrow\R$. If $h\in C(M,\R)$ then $\mu(h)=\int_Mhd\mu$. Let us denote by $\rho$ the metric on $\SM(f)$ which defines the weak-* topology as follows. Let $\phi_1, \phi_2,\ldots$ be a countable dense subset of the unit ball in $C(M,\R)$ and set
$$\rho(\mu,\nu)= \sum_{i\geq1} \frac{1}{2^i} |\mu(\phi_i) - \nu(\phi_i)|.$$

Given a periodic ergodic measure $\mu_p \in \SM_e(f)$, we denote by $\chi^{+}(p,f)$ and $\chi^{-}(p,f)$ the smallest positive Lyapunov exponent and the biggest negative Lyapunov exponent of $\mu_p$, respectively. Then we define $\chi(p,f)=\min\{\chi^+(p,f),\, -\chi^-(p,f)\}$.% the positive characteristic exponent of $\mu$, which is defined by Osedelets' theorem: for $\mu$-almost every point $x$, $$\lim_{n\rt \y}\frac{1}{n} \log |Df^n_x|= \chi(\mu),$$\margem{OLhar a notacao de Derivada!!!}

%If $p$ is a periodic point, we denote by $\chi(p)$ the positive characteristic exponent of the periodic measure $\mu_p$.
%\margem{este $\chi(f)$ nao e assim nao! e talvez seja melhor ele aparecer depois!}
%We also define  $$\chi(f)= \inf\{\chi(p) : p \in Per_h(f) \,\ and \,\ \tau(p)=\tau(f)\}.$$

%\subsection{Partitions and Entropy}

%\subsection{Hausdorff Topology}

%Let $\rho_H$ denote the \textit{Hausdorff metric} on the collection of compact subsets of $\SM(f)$. DEFINIR!

%This paper is organized as follows: %explicar como vai ser o paper
\section{The property $\SS_{n,p}$}\label{SNP}

In this section, we define and study the $\SS_{n,p}$ property. This property is in the spirit of \cite{D-N}, %in order to use  the entropy structure, given by Boyle and Downarowicz \cite{BD}, 
in order to find diffeomorphisms that has no symbolic extensions.

\begin{definition}
Given a positive integer $n$, we say that a diffeomorphism $f$ satisfies property $\SS_{n,p}$ if $p$ is a hyperbolic periodic point of $f$, and for any $\tilde{p} \in Per_h^n(f)$ related to $p$ there is a zero dimensional periodic hyperbolic basic set $\La(\tilde{p}, n)\subset H(p,f)$ for $f$ with the same index that $p$,
 such that the following happens:
\begin{itemize}
 %\item[a)] there is a zero dimensional periodic hyperbolic basic set $\La(\tilde{p}, n)\subset H(p,f)$ for $f$ with the same index that $p$,
 %such that $$\La(\tilde{p}, n) \cap \partial\al_n = \es,$$
%\item[b)] $\La(\tilde{p}, n)$ is subordinate to $\al_n,$
\item[a)] there is $\nu \in \SM_e(\La(\tilde{p}, n))$ such that $$h_\nu(f) >\chi(\tilde{p},f)- \frac{1}{n},$$
\item[b)] for every $\mu \in \SM_e(\La(\tilde{p}, n))$, we have $$\rho(\mu,\mu_{\tilde{p}}) <\frac{1}{n}.$$
%\item[e)] for every $\mu \in \SM_e(\La(\tilde{p}, n))$ we have $$\chi(\mu)>\chi(\tilde{p})-\frac{1}{n}.$$
\item[c)] for every hyperbolic periodic point $q \in \La(\tilde{p}, n)$,  we have $$\chi(q,f)>\chi(\tilde{p},f)-\frac{1}{n}.$$
\end{itemize}
\end{definition}

%From now on we will consider $f \in \mundo$, $H(p)$ a $f$-homoclinic class associated with the $f$-%hyperbolic periodic point $p$, and $\SU(f)$ any small neighborhood of $f$ such that for all $g\in\SU(f)$ there is a continuation of $p$, that we will denote by $p_g$

The following result concern about the abundance of diffeomorphisms satisfying property $S_{n,p}$ near diffeomorphism with homoclinic classes admitting no dominated splittings.

\begin{proposition}
Let $f$ be a Kupka-Smale generic diffeomorphism with a hyperbolic periodic point $p$ of index $i$. If $H(p,f)$ is a non-trivial homoclinic class admitting no $i-$dominated splitting,  then for any neighborhood $\mathcal{U}$ of $f$ and any positive integer $n$, there exists an open subset $\mathcal{V}\subset\mathcal{U}$ such that every $g\in \mathcal{V}$ satisfies property $\SS_{n,p(g)}$.
\label{mainprop}\end{proposition}

The idea to prove this Proposition is to produce many nice horseshoes, as done by Downarowicz and Newhouse in \cite{D-N}. However, in their context, there is an abundance of homoclinic tangencies to produce such horseshoes. In our context we will use Lemma \ref{mainlemma}, which is a key and technical lemma,  to overcome the lack of such abundance in general.

\begin{proof}
We can suppose that every periodic orbit of $Per_h^n(f)$ and the orbit of $p$ has an analytic continuation on $\SU$. Moreover, by semicontinuity arguments, there exists $k$ such that for every $g\in \SU$ we have
$$\#\{q\in Per_h^n(g);\textrm{ homoclinically related with $p$}\}=k.$$We denote the elements of this set for $f$ by
$\{p_1,\dots,p_k\}$.

%Given $f\in \mathcal{R}$ be a $C^1$ diffeomorphism, $p$ a hyperbolic periodic point of $f$ and  $H(p,f)$ as in the
%hypothesis, let $\mathcal{U}$ be a small neighborhood of $f$ and $n$ a positive integer fixed arbitrary.

%Since $\mathcal{R}$ is the Kupka-Smale residual subset, we can suppose $\mathcal{U}\subset \mathcal{A}_n$, and moreover such that the cardinality of $Per_h^n(g)\cap H(p(g),g)$ is finite and constant for diffeomorphisms inside $g\in\mathcal{U}$.  Hence, let $p_1,\, \ldots,\, p_k$  the hyperbolic periodic points of $f$, homoclinic related to $p$ with period smaller or equal than $n$.

Since $H(p,f)=H(p_1,f)$ admits no $i-$dominated splitting, Gourmelon's result \cite{G} implies that, after some perturbation, we can suppose that $f$ exhibits a homoclinic tangency for $p_1$, i.e., there exists a non
transversal intersection between $W^s(O(p_1),f)$ and $W^u(O(p_1),f)$.

%%%%%%%%% Ver isto depois para o caso geral %%%%%%%%%%%
%Let $f\in \mathcal{R}$ be a $C^1$ diffeomorphism and $p$ a hyperbolic periodic point such that $H(p,f)$ is  %not a partial hyperbolic set of
%$f$.  Then, by Theorem Blá (ESTE TEOREMA É O TEOREMA DE CROVISIER-%SAMBARINO-YANG) after some perturbation, we can suppose there exists a
%hyperbolic periodic point %%$p_1\in H(p,f)$ with index $i$ such that $W^s(o(p_1),f)$ has a non transversal intersection with $W^u(o(p_1),f)
%$. Note, the index of $p_1$ could be different of index of $p$.

%\begin{lemma}
%Let $f$ be a $C^1$ diffeomorphism  with a homoclinic class $H(p,f)$. Suppose  there exist a hyperbolic periodic point $p_1$ exhibiting a homoclinic tangency, and  a hyperbolic periodic point $p_2$ of $f$ related to $p_1$.
%Given any small neighborhood $\mathcal{V}$
%of $f$,  if $n$ is large enough, there exist a diffeomorphism
%$g\in \mathcal{V}$, a hyperbolic set $\Lambda(p_1(g),n)\subset H(p(g),g)$ of $g$, such that the items of  property $S_{n,p(g)}$ are true for $p_1(g)$, and moreover $g$ exhibits a homoclinic tangency for $p_2(g)$. %Here $p_1(g)$ and $p_2(g)$ are the continuations of
%the hyperbolic periodic points $p_1$ and $p_2$ for $g$.
%\label{mainlemma}\end{lemma}

%Now we state two technical lemmas.

Now we state a technical lemma.

\begin{lemma}
If $n$ is large enough, there exist a diffeomorphism $g\in \SU$ and a small neighborhood $\mathcal{V}\subset\SU$ of $g$ such that for every $h\in \mathcal{V}$ the items of property $S_{n,p(h)}$ holds for $p_1(h)$, and moreover $g$ exhibits a homoclinic tangency for $p_2(g)$. %Moreover, there exists an open neighborhood $\SW\subset \SU$ of $g$ such that for every $h\in \SW$ the items of the property $S_{n,p(h)}$ holds for $p_1(h)$.
\label{mainlemma}\end{lemma}

We postpone the proof of this lemma and finish the proof of the Proposition.

We consider now $g_1$ and the neighborhood $\mathcal{V}_1$ of $g_1$ given by Lemma  \ref{mainlemma}. Now, since $g_1$ exhibits a homoclinic tangency for $p_2(g_1)$, and $p_3(g_1)$ is homoclinically related to $p_2(g_1)$, we can use Lemma \ref{mainlemma} again to obtain a diffeomorphism $g_2$ and a neighborhood $\mathcal{V}_2\subset \mathcal{V}_1$ of $g_2$ such  that for every diffeomorphism $h\in \mathcal{V}_2$ the items of property $S_{n,p(h)}$ holds for $p_2(h)$, and moreover $g_2$ exhibits a homoclinic tangency for $p_3(g)$.

Now, we repeat the process finitely many times, to obtain a diffeomorphism $g=g_k$ and a neighborhood $\mathcal{V}=\mathcal{V}_k\subset \mathcal{V}_{k-1}\ldots \subset \mathcal{V}_1\subset \SU$ of $g$ such that the items of property $S_{n,p(h)}$  holds for $p_i(h)$ with $i=1,\dots, k$ for every $h\in \mathcal{V}$. And then, by choice of $\SU$, every diffeomorphism $h\in \mathcal{V}$ satisfy property $S_{n,p(h)}$.

%Using Lemma \ref{mainlemma} a hyperbolic set $\Lambda(p_1(g),n)$ is created satisfying the items of property $S_{n,p(g)}$. Moreover, this set persists and still satisfies the items. However, since $g$ still have an homoclinic tangency for $p_2(g)$, we can use Lemma \ref{mainlemma} again to obtain $h\in \SW$ such that the items of property $S_{n,p(h)}$ holds for $p_2(h)$ and $p_3(h)$ has a homoclinic tangency. By the persistence mentioned before, the items of property $S_{n,p(h)}$ also holds for $p_1(h)$. Now, we repeat the process finitely many times, to obtain a diffeomorphism $\widetilde{f}\in \SU$ such that the items of the $S_{n,p}$ property  holds for $p_i(\widetilde{f})$ with $i=1,\dots, k$. The persistence mentioned before  gives the open set $\SV$.

%First, namely, the hyperbolic set built in the proof of the previous lemma is such that the items of the property $S_{n,p}$ is open. See Remark \ref{open}.  Now, since $Per_h^n(g)\cap H(p,g)$ is constant for diffeomorphisms $g\in \mathcal{U}$,  we can apply the above lemma finitely many times, to find a diffeomorphism $g\in \mathcal{U}$ satisfying the property $S_{n,p(g)}$.  And then, by choice of $\mathcal{U}$ and by Remark \ref{open} once more, $g$ satisfies robustly property $S_{n,p}$.  Which concludes the proposition.

\end{proof}

%\begin{lemma}
%Let $f$ be a $C^1$ diffeomorphism exhibiting a homoclinic tangency for a hyperbolic periodic point $p_1$. Given a neighborhood $\mathcal{U}$
%of $f$, a hyperbolic periodic point $p_2$ of $f$ homoclinic related to $p_1$ there exists $n_0$ with the following property. If  $n\geq n_0$, there exist a diffeomorphism
%$g\in \mathcal{U}$, a hyperbolic set $\Lambda(p_1(g),n)\subset H(p_1(g),g)$ of $g$ satisfying properties (1), (2), (3), (4) and
%(5) of property $S_{n,m,p}$????, and moreover $g$ exhibits a homoclinic tangency for $p_2(g)$. Here $p_1(g)$ and $p_2(g)$ are the continuations of
%the hyperbolic periodic points $p_1$ and $p_2$ for $g$.
%\label{mainlemma}\end{lemma}

%AQUI AGORA, USA O LEMA PARA TERMINAR A PROVA DA PROPOSIÇÃO. DE QUE MANEIRA: VOCÊ PODE SUPOR QUE A VIZINHANÇA $U$ SEJA TAL QUE ELA ESTEJA CONTIDA
% NO CONJUNTO $\mathcal{A}_n$, ONDE ESTE É O CONJUTO DADO POR KUPKA SMALE, QUE TODOS OS PONTOS PERIÓDICOS DE PERÍODO MENOR OU IGUAL A $n$ SÃO
%HIPERBÓLICOS. ASSIM, SÃO FINITOS. AH, A VIZINHANÇA TAMBÉM TEM Q SER TOMADA ONDE A CLASSE HOMOCLÍNICA TEM CONTINUIDADE E TAL, PARA VC GARANTIR
%QUE A QUANTIDADE DE PONTOS PERIÓDICOS NA CLASSE DE PERÍODO MENOR OU IGUAL A $n$ É FINITA E CONSTANTE NESTA VIZINHANÇA. PRONTO, DESDE QUE AS
%PROPRIEDADES SÃO ABERTAS COMO JÁ FOI VISTO ACIMA, VC REPETE O LEMA PARA ESTES FINITOS PONTOS E GANHAMOS.

\subsection{Proof of Lemma \ref{mainlemma}}
%\margem{ainda não revisei a prova do lema então a notação está diferente das anteriores. Ei, nao tem que citar a zeze nao? pq o metodo da criacao da ferradura eh parecido}

First of all, we observe that many times in this proof we use expressions like ``by some
perturbation", or ``we can perturb $f$", to say we can take a diffeomorphism arbitrary close to $f$. Sometimes, in order to not complicate the notation we use the same letter to denote the new  diffeomorphism. Also, when we say ``by a local perturbation" we mean
that we can perform a perturbation of $f$ keeping the new diffeomorphism equal to $f$ outside some small open set.

%First of all, we observe that many times in this proof we use expressions like ``by some
%perturbation", or ``we can perturb $f$", which means that arbitrary close to $f$  there exists a diffeomorphism with a desire property. Sometimes, in order to not complicate the notation we will use the same notations for the new  diffeomorphism, we hope that this does not cause confusions. Also, when we say ``by a local perturbation" we mean
%that we can make one perturbation of $f$ such that outside some small open set, the new diffeomorphism remains equal  to $f$.

Let $q$ be a point of homoclinic tangency of $p_1$, and $V$ be a small neighborhood of $O(p_1)$ such that
$f^{-1}(q)$ is not in $V$. Shrinking $V$, if necessary, we can suppose $f^{\tau(p_1,f)}=Df^{\tau(p_1,f)}$ (in local coordinates on $V$) after
a perturbation (see Franks' lemma \cite{F}). We remark that after this perturbation the homoclinic tangency could disappear. Nevertheless, since $f^{-1}(q)$ is not in $V$, using the continuity of compact parts of unstable and stable manifolds of $p_1$,  by a local perturbation in some neighborhood of $f^{-1}(q)$  we can
recover the homoclinic tangency.

%Now, we replace $q$ for the point of homoclinic tangency such that $q\in V$ and $f^{-1}(q)\not\in V$. After that, we take a neighborhood $U$ of
%$q$ such that $f^{-1}(U)\cap V=\emptyset$.  We denote by $D$ the connected component of $W^u(p,f)\cap U$ that contains $q$.

Up to take another point of the orbit of $q$, we can suppose that $q\in V$ and $f^{-1}(q)\not\in V$. So, we can take a neighborhood $U$ of
$q$ such that $f^{-1}(U)\cap V=\emptyset$.  We denote by $D$ the connected component of $W^u(p_1,f)\cap U$ that contains $q$.

Now, we look to $U$ in some local coordinates with the  splitting $T_qD\oplus T_qD^{\bot}$, and such that $q=0$  in these coordinates.
Since $D\subset W^u(p_1,f)$ we have that $D$ is a graph of a $C^1$  map $r:T_{q}D\rightarrow T_{q}
D^{\bot}$, i.e. $D=(x,r(x))$.  Moreover, $Dr(q)$ is close to zero.  Hence, the
diffeomorphism $\phi(x,y)=(x,y-r(x))$ is $C^1$ close to identity in a small  neighborhood of $q$. In particular, there exists  a diffeomorphism  $h$,  $C^1$-close to identity, such that $h=\phi$ in some small neighborhood of $q$, and $h=Id$ for points far away from $q$. Thus, $f_1:=h\circ f$ is a $C^1$ local perturbation of $f$ such that
$T_qD\cap U\subset W^u(p_1,f_1)$. Since $f^{-1}(U)\cap V=\emptyset$, we have that $f_1=f$ in $V$, as a consequence $f_1|V$ is still linear, and
$W^s_{loc}(p_1,f_1)$ remains unchanged in $U$.  Since $q$ is a non transversal homoclinic
point we have that $T_q D\cap E^s(p_1,f)$ is a non trivial subspace. Actually, we can assume that
$T_q D\cap E^s(p_1,f)$ is an one-dimensional subspace, after some local perturbation if necessary. Thus,  $f_1$ exhibits an interval of homoclinic tangencies containing $q$.
%Now, since $T_q D\cap E^s(p_1,f)$ is a non trivial subespace ($q$ is a non transversal homoclinic
%point), the diffeomorphism $f_1$ contains an interval of homoclinic tangencies containing $q$. We remark that we could have assumed that
%$T_q D\cap E^s(p_1,f)$ is an one-dimensional subspace, after some local perturbation.

Let $I$ be this interval of homoclinic tangencies. Replacing the local coordinates in $U$, if necessary,  we can suppose
that $\{(x_1,0,...,0), \, -3a\leq x_1\leq 3a\} \subset I$, for some $a>0$ small enough. % and usual coordinates of
%$\R^{d}$.

Let $N$ be a large positive integer. Taking $I$ smaller, if necessary, we can construct a diffeomorphism $\Theta:M\rightarrow M$, such that $\Theta=Id$ in  $B(0,2a)^c$ and
$$
\Theta(x,y)=\left(x_1,...,x_s,\,y_1+A\cos\frac{\pi x_1 N}{2a},\, y_2,...,y_u\right), \text{ for } (x,y)\in B(0,a)\subset U,
$$
for $A=\displaystyle\frac{2Ka\delta }{\pi N}$, where $K$ is a constant which depends only on the local coordinates over $U$ and $\delta>0$ is so
small as we want. %The argument is the same as the one that we used to find $h$, using bump functions.

Hence, taking   $g=\Theta\circ f_1$, we have that $g$ is $\delta-C^1$ close to $f_1$ and moreover $g=f_1$ in the complement of
$f_1^{-1}(B(q,2a))$.  Note that $g$ depends on $N$ but to not complicate the notation we denote this diffeomorphism by $g$,  independent of
$N$.

\begin{remark}
The most important properties of this new diffeomorphim  is that $g$ has $N$ transversal homoclinic points for $p_1$ inside $U$, but  $g$ still
has an interval of homoclinic tangency inside $U$, in fact there are two intervals of homoclinic tangency in $U$: one inside
$\{(x_1,0,...,0), \, -3a\leq x_1\leq -2a\}$ and other inside $\{(x_1,0,...,0), \, 2a\leq x_1\leq 3a\} $, in local coordinates.
\label{tangency}\end{remark}

To simplify notation we assume $p_1$ is a fixed point, being similar the general case.

We remark that $g|V$ is still linear in local coordinates, since $f$ is equal $g$ in $V$. Let $D_t=D^s\times D^u_t$ be a small rectangle, with $D^s=W^s_{loc}(p_1,g)\cap U$, and $D^u_t$  a small disk in
$\{(0,\ldots,0,y_1,\ldots,y_n),\, y_i\in\R^+ \text{ and } |y_i|<A/4\}$, such that  $t$ is the smallest positive integer such that $g^t(D_t)$ is a
disk $A/4-C^1$ close to the connected component of $W^u(p_1,g)\cap U$ containing the $N$ transversal homoclinic points built before.  We remark
 that $t$ depends on $N$, and $t\to\infty$ when $N\to \infty$.

Observe that $A$ is small if $N$ is large, and by choice of $D_t$, we have that $g(D_t)\cap D_t$ has $N$ disjoint connected components. Moreover, %taking $N$ larger (which implies $A$ smaller), if
%necessary,
we have that the maximal invariant set in $D_t$ for $g^t$
$$
\tilde{\Lambda}(p_1,N)=\bigcap_{j\in\Z} g^{tj}(D_t)
$$
is a hyperbolic set inside $H(p_1,g)$. %This is the so called Newhouse's snake.

Let $\Lambda(p_1,N)={\displaystyle\bigcup_{0\leq j\leq t} g^j(\tilde{\Lambda}(p_1,N))}$ be the hyperbolic periodic set of $g$ induced by
$\tilde{\Lambda}(p_1,N)$.
Since $g|\tilde{\Lambda}(p_1,N)$ is conjugated with the full shift of
$N$ symbols,  we have that $h(g|\Lambda(p_1,N))=\displaystyle\frac{1}{t}\log N$,

We recall that $g|V$ is linear. So, if $m$ is the largest positive integer such that $g^j(x)\in V$ for $0\leq j\leq m$, there
exist constants $K_1$ and $K_2$ depending on the local coordinate on $V$  such that
\begin{equation}
K_1\|Dg(p_1)^{m}|E^u\|^{-1}\leq d(x,W^s_{loc}(p_1,g))\leq K_2 \|Dg(p_1)^{-m}|E^u\|, \label{desigualdade 1}\end{equation} for $x\in V$.
Analogously, if $m$ is the largest positive integer such that $g^{-j}(x)\in V$ for $0\leq j\leq m$, then  there exist constants  $K_3$ and $K_4$ such that
\begin{equation}
K_3 \|Dg(p_1)^{-m}|E^s\|^{-1} \leq d(x,W^u_{loc}(p_1,g))\leq K_4 \|Dg(p_1)^{m}|E^s\|. \label{desigualdade 2}\end{equation}

Another consequence of $g|V$ be linear and the choice of $t$ is the following result, which also appears in
\cite{C-T}.

\begin{lemma}
[Lemma 4.2 of \cite{C-T}] \label{afirma}
For $A$ and $t$ defined as before, there exists  a positive integer $K_5$, which is %and $K_2$
independent of $A$, such that $$A<K_5\max \{ \|Dg(p_1)^{-t}|E^u\|,\,
\|Dg(p_1)^{t}|E^s\|\}.$$\end{lemma}

Let $n$ be a large positive integer. Since $A=\displaystyle\frac{2Ka\delta}{\pi N}$,  using Lemma \ref{afirma} and recalling that $N\to \infty$ implies $t\to \infty$, we can select a large positive integer $N$, such that
$$
\displaystyle\frac{1}{t}\log N> \min\left\{\frac{1}{t}\log\|Dg(p_1)^{-t}|E^u\|^{-1},\, \frac{1}{t}\log\|Dg(p_1)^{t}|E^s\|^{-1} \right\} -\frac{1}{2n}.
$$
%It is clear that $t$ goes to infinity when $N$ goes to infinity.
But, when $t$ goes to infinity the above minimum converges to $\chi(p_1,g)$,
by definition. Therefore, there exists a large positive integer $N_1$ such that
$$
\displaystyle\frac{1}{t}\log N_1> \chi(p_1,g)-\frac{1}{n}.
$$
So, it is possible to find a $C^1-$perturbation $g$ of $f$ such that
$$
h(g|\Lambda(p_1,N_1))> \chi(p_1,g)-\frac{1}{n}.
$$

Now, by the variational principle there exists an
ergodic measure $\mu_N\in\mathcal{M}(\Lambda(p_1,N))$ such that
\begin{equation}
h_{\mu_N}(g)> \chi (p_1,g)-\frac{1}{n}, \text{ for } N\geq N_1.\label{item b}\end{equation}

Observe that the orbit of points in the hyperbolic set $\Lambda(g,N)$, when $N$ is large enough, stay almost all the time inside the neighborhood
$V$ of $p_1$, which one could be assumed so small as we wanted. Hence,  there exists a positive integer
$N_2$ such
that if $\mu\in \mathcal{M}(f|\Lambda(g,N))$ is ergodic then $\rho(\mu,\mu_{p_1})<1/n$, for every $N\geq N_2$.

Finally, we find $N_3$ in order to obtain property (e) of $S_{n,p(g)}$ for $\Lambda(p_1,N)$ with $N\geq N_3$.

We define
$$
V_k^{u}=V\cap g(V)\cap...\cap g^{k}(V), \text{ and }
$$
$$
V_k^{s}=V\cap g^{-1}(V)\cap...\cap g^{-k}(V).
$$

Given vectors $v,w\in\R^{2n}$ and subspaces  $E,F\subset\R^{2n}$ we define
$$
ang(v,w):=\left|\tan\left[\arccos\left(\frac{<v,w>}{\|v\|\|w\|}\right)\right]\right|,$$

$$
ang(v,E)=\min_{w\in E,\, |w|=1}\, ang(v,w)\quad \text{and}\quad ang(E,F)=\min_{w\in E,\, |w|=1} \, ang(w,F).
$$

The following lemma, is also  a straightforward consequence of $g|V$ be linear, as in Lemma 4.4 in \cite{C-T}.

\begin{lemma}%[Lemma 4.4 in \cite{C-T}]
With above definitions, there exists positive constants  $K_6$ and $K_7$, such that
\begin{itemize}

\item[1)] if $z\in V^{u}_k$, $v\in \R^{2n}\backslash E_{p_1}^u$ and $ang(g^{-k}(v),E_{p_1}^u)\geq
1$, then
$$
K_6 \|Dg_{p_1}^{k}|E^s\|^{-1}|v| \,\min \{ang(v,E_{p_1}^{u}),\, 1\}\leq |Dg^{-k}(z)(v)|\leq K_7 \|Dg_{p_1}^{-k}|E^s\||v|
$$

\item[2)] if $z\in V^{s}_k$, $v\in \R^{2n}\backslash E_{p_1}^s$ and $ang(g^{k}(v),E_{p_1}^s)\geq
1$, then
$$
K_6 \|Dg_{p_1}^{-k}|E^u\|^{-1}|v| \,\min \{ang(v,E_{p_1}^{s}),\, 1\}\leq |Dg^{k}(z)(v)|\leq K_7 \|Dg_{p_1}^k|E^u\||v|
$$
\end{itemize}

\label{lema}\end{lemma}

Now,  since $\Lambda(p_1,N)=\cup_{i=0}^{t-1} g^i(\tilde{\Lambda}(p_1,N))$
%is a hyperbolic set for $g$,
with $\tilde{\Lambda}(p_1,N)\subset
V$, then we can take positive integers $k$ and $T$  such that $t=k+T$, and  $g^{i}(\tilde{\Lambda}(p_1,N))\subset V$ for $0\leq i\leq k$. Moreover, by construction of $\tilde{\Lambda}(p_1,N)$ this $T$ can be taken independent of $N$. Hence, provided $t$ goes to infinity when $N$
goes to infinity, we have that $k$ also goes to infinity.  Now,  %$\tilde{\Lambda}(p_1,N))$
we know that the hyperbolic decomposition $T_{\tilde{\Lambda}(p_1,N)}M=\tilde{E}^s\oplus \tilde{E}^u$ of the hyperbolic set $\tilde{\Lambda}(p_1,N)$ is such that $\tilde{E}^s(g^{-k}(z_1))$ and $\tilde{E}^u(g^{k}(z_2))$ are
close to $E^s_{p_1}$ and $E^u_{p_1}$, respectively, for every $z_1\in g^{k}(\tilde{\Lambda}(p_1,N))$ and $z_2\in\tilde{\Lambda}(p_1,N)$. In
particular, 
$$ang(Dg^{-k}(z_1)(v),E^u_{p_1})>1\textrm{ for }v\in \tilde{E}^s(z_1)\textrm{ and }$$ $$ang(Dg^{k}(z_2)(v),E^s_{p_1})>1\textrm{ for }v\in \tilde{E}^u(z_2).$$
Moreover, and the most important argument in this case, is that  although $ang(v,E_{p_1}^{u})$ for $v\in \tilde{E}^s(z_1)$, and $ang(v,E_{p_1}^{s})$
for $v\in \tilde{E}^u(z_2)$ are a very small constant, independent of $N$, we have ensured that 
$$ang(Dg^{-k}(z_1)(v),E^u_{p_1})>1\textrm{ and }ang(Dg^{k}(z_2)(v),E^s_{p_1})>1.$$

So, using these informations and Lemma \ref{lema} we can find  constants $K_6$ and $K_7$, such that for every $z\in \tilde{\Lambda}(p_1,N)$,
 $r=l(k+T)$ and for every $l\in\N$:
 \begin{itemize}
\item[1)] if $v\in \tilde{E}^s(z)$ then
\begin{equation*}
 |Dg^{-r}(z)(v)| \geq (C_1\, K_6)^l \, \|Dg_{p_1}^{k}|E^s\|^{-l}|v|  \label{equa 2}\end{equation*}

\item[2)] if $v\in \tilde{E}^u(z)$ then
\begin{equation*}
|Dg^r(z)(v)|\geq (C_1\, K_7)^l \, \|Dg_{p_1}^{-k}|E^u\|^{-l}|v|, \label{equa 2}\end{equation*}
\end{itemize}
where
$$C_1=\inf_{z\in V\backslash g^{-1}(V), \; |v|=1} \|Dg^{T}(z)(v)\|.$$

Therefore,  for $N$ large enough, all points in $\tilde{\Lambda}(p,N)$ have Lyapunov exponents with absolute values bigger than
$\chi(p,g)-1/n$. In particular, we can choose $N_3$, in order to get $k>>T$, such that for any periodic point $\tilde{q}\in \Lambda(p_1,N)$, with $N>N_3$, we have
$$\chi(\tilde{q})> \chi(p_1,g)-\frac{1}{n}.$$
Hence, if we take  $\Lambda(p_1,n)=\Lambda(p_1,N)$ for $N=\max\{N_1,\, N_2,\, N_3\}$, the items of property $S_{n,p(g)}$ are satisfied for the
perturbation $g$ of $f$ and the hyperbolic periodic point $p_1$ of $g$.

\begin{remark}
By the construction of $\Lambda(p_1,n)$, observe that every item in the property $S_{n,p}$ is robust. That is, if $\tilde{g}$ is close to $g$, then the continuation of $\Lambda(p_1,n)$ for $\tilde{g}$  is such that all the items of property $S_{n,p(\tilde{g})}$ is still true for $p_1(\tilde{g})$. This is because item (a) is a robust property; item (b) is a consequence of the continuation of $\Lambda(p_1,n)$ be also inside $V$ and  item (c) still is true for continuations of $\Lambda(p_1,n)$ by Lemma \ref{lema}.
\label{open}\end{remark}

Hence, by the previous remark there exists a neighborhood $\mathcal{V}\subset \mathcal{U}$ of $g$, such that every diffeomorphism $h\in\mathcal{V}$ satisfies the items of property $S_{n,p(h)}$ for $p_1(h)$.

Now,  since the diffeomorphism $g$ belongs to $\mathcal{U}$, we know that the hyperbolic periodic point $p_2(g)$ still is
homoclinic related with $p_1(g)$. Also, by Remark \ref{tangency}, $g$ still exhibits a homoclinic tangency for $p_1(g)$. Now, by a perturbation using Franks Lemma, we can find a transversal homoclinic point to $p_1(g)$, such that the angle between $W^s(p_1(g),g)$ and $W^u(p_1(g),g)$ is so small as we want.  Hence, since $p_1(g)$ and $p_2(g)$ are related, there exists a transversal homoclinic point for $p_2(g)$ such that the angle between $W^s(p_2(g),g)$ and $W^u(p_2(g),g)$ is so small, too. Finally, using Franks Lemma once more, we can perturb $g$ such that this transversal homoclinic point become a homoclinic tangency. Since this perturbation can be find in $\mathcal{V}$, we finish the proof.

%Finally Remark \ref{open}  gives the existence of the open set $\SV$.

\section{A local version of the main theorem}\label{localversion}

First we recall some knowns residual subsets. We denote by $\mathcal{R}_1\subset\diff^1(M)$ the residual subset given by \cite{C-M-P}, such that for every diffeomorphism $g\in \mathcal{R}_1$ two
homoclinic classes are either disjoint or coincide. By $\mathcal{R}_2\subset \diff^1(M)$ the residual subset given by \cite{A-B-C-D-W},
such that for every diffeomorphism $g\in \mathcal{R}_2$, every homoclinic class having a hyperbolic periodic point with index $i$ and a hyperbolic periodic point with index $j$, with $i<j$, has  a dense set of
hyperbolic periodic points with index $k$ for every $i\leq k\leq j$. And by $KS$ the residual subset of Kupka-Smale diffeomorphisms. %and by $\mathcal{R}_3$ the residual subset given by Proposition \ref{prop4.4}. 
Hence, we define $\mathcal{R}_4=\mathcal{R}_1\cap\mathcal{R}_2%\cap\mathcal{R}_3
\cap KS$.

%Before we prove this proposition, let us use it to prove one local version of the main theorem.

\begin{proposition}\label{local}
Let $f \in \mathcal{R}_3$, $p$ be a hyperbolic periodic point of $f$. If  $\SU(f)\subset\mundo$ is a small enough neighborhood of $f$, there is a residual subset
$\mathcal{R} \subset \SU(f)$ such that every $g \in \mathcal{R}$ satisfies only one of the following statements:
\begin{itemize}
 \item[(i)] $H(p_g,g)$ has a good decomposition;
\item[(ii)] $g|_{H(p_g,g)}$ has no symbolic extensions.
\end{itemize}
%Here, $H(p_g)$ denotes the continuation of $H(p,f)$.
\end{proposition}

\textbf{Proof:}

%To find such residual subset we use first some known residual subsets. We denote by $\mathcal{R}_1\subset\diff^1(M)$ the residual subset given by \cite{C-M-P}, such that for every diffeomorphism $g\in \mathcal{R}_1$ two
%homoclinic classes are either disjoint or coincide. By $\mathcal{R}_2\subset \diff^1(M)$ the residual subset given by \cite{A-B-C-D-W},
%such that for every diffeomophism $g\in \mathcal{R}_2$, every homoclinic class having a hyperbolic periodic point with index $i$ and a hyperbolic periodic point with index $j$, with $i<j$, has  a dense set of
%hyperbolic periodic points with index $k$ for every $i\leq k\leq j$.  And, by $\mathcal{R}_3$ the residual subset given by Proposition \ref{prop4.4}.

% and finally we denote by $\mathcal{R}_4$ the residual subset in $\mathcal{U}(f)$ given by Lemma \ref%{residual}.

Since $f\in \mathcal{R}_3$ if $\mathcal{U}(f)$ is small enough then there exist $i$ and $j$, and hyperbolic periodic points  $p_i, \, p_{i+1},\ldots, \, p_{j}$ of $f$ with $ind\, p_k=k$, for $i\leq k\leq j$, such that:
\begin{itemize}
\item[-] $H(p,f)=H(p_i,f)=H(p_{i+1},f)=\ldots=H(p_j,f)$,

\item[-] for every hyperbolic periodic point $q\in H(p,f)$ we have $i\leq ind\, q\leq j$.
\end{itemize}

By \cite[Lemma 4.2, pg.20]{A-B-C-D-W}, there exists an open and dense subset of $\mathcal{U}(f)$ over $\mathcal{R}_3$ such that for every $g$ in this subset, we still have that:
\begin{itemize}
\item[-] $H(p(g),g)=H(p_i(g),g)=H(p_{i+1}(g),g)=\ldots=H(p_j(g),g)$,

\item[-] for every hyperbolic periodic point $q\in H(p(g),g)$ we have $i\leq ind\, q\leq j$.
\end{itemize}

%Now, we will find some residual subset in $\mathcal{U}(f)$.
Now, for any positive integer $n$ and any $i\leq k\leq j$, we define $\mathcal{B}_{n,p_k}\subset \mathcal{U}(f)$ as the subset of diffeomorphisms that robustly satisfies property $S_{n,p_k}$, i.e.,  $g\in\mathcal{B}_{n,p_k}$ if there is a small neighborhood of $g$ where every diffeomorphism $h$ satisfy property $S_{n,p_k(h)}$.

\begin{lemma}\label{residual}
There is a residual subset of $\mathcal{R}_4\subset \mathcal{U}(f)$, such that for any  positive integer $n$ and any $i\leq k\leq j$, if $g \in \mathcal{R}_4$ and there is a sequence of diffeomorphisms $\{g_m\}\in \mathcal{B}_{n,p_k}$
that converges to $g$,  then $g$ satisfies property $\SS_{n,p_k(g)}$.
\end{lemma}
\proof
Let us define $\mathcal{V}_{n,k}=\mathcal{B}_{n,p_k}\cup \overline{\mathcal{B}_{n,p_k}}^c$ be an open and dense subset in $\mathcal{U}(f)$, for every positive integer $n$ and every $i\leq k\leq j$. Then,  $\mathcal{R}_4=\cap_{n\geq 0}\cap_{i\leq k\leq j} \mathcal{V}_{n,k}$ is a residual subset in $\mathcal{U}(f)$. To finish the proof, let $g\in \mathcal{R}_4$. Given a positive integer $n$ and $i\leq k\leq j$, if there exists diffeomorphisms $g_m\in \mathcal{B}_{n,p_k}$ converging to $g$, then $g\not\in\overline{\mathcal{B}_{n,p_k}}^c$. Therefore, since $g\in \mathcal{V}_{n,k}$ we have that $g\in\mathcal{B}_{n,p_k}$ and then satisfies property $S_{n,p_k(g)}$.

$\hfill\square$

Using Lemma \ref{residual}, we define $\mathcal{R}=\mathcal{R}_3\cap \mathcal{R}_4$, which is a residual subset in $\mathcal{U}(f)$. Now, we will verify that a diffeomorphism in this residual subset satisfies one of the two properties claimed in the proposition which finishes the proof.

For this we will use the following result of Burguet.  %the proof could be simplified. %and thus just for sake of completeness we will give a sketch here. 

\begin{proposition} [Corollary 1 of \cite{B1}]
Let $f:M\to M$ be a dynamical system admitting a symbolic extension. Then the entropy function $h: \mathcal{M}(f)\rightarrow \R$ is a difference of nonnegative upper semicontinuous functions. In particular the entropy function $h$ restrict to any compact set of measures has a large set of continuity points. 
\label{propburguet}\end{proposition}

Let $g\in \mathcal{R}$, by choice of $\mathcal{R}$, $i$ and $j$ are the two extreme indices in $H(p(g),g)$ and
$H(p(g),g)=H(p_i(g),g)=H(p_{i+1}(g),g)=\ldots=H(p_j(g),g)$.

Suppose $H(p(g),g)$ admits no good decomposition. Hence,  there is some $i\leq k\leq j$ such that $H(p(g),g)=H(p_k(g),g)$ admits no $k-$dominated splitting.  %Since $g$ belongs to a residual subset, we could have assumed that $g$ is a Kupka-Smale generic diffeomorphism.

By Proposition \ref{mainprop}, for every $n>0$, we can find a sequence of diffeomorphisms $\{g_{n,m}\}_{m\in \N}$  converging to $g$,  such that each $g_{n,m}\in \mathcal{B}_{n,p_k}$.  Therefore, by Lemma \ref{residual}, $g$ satisfies property $S_{n,p_k(g)}$ for every $n>0$, since $g\in \mathcal{R}_4$.

%Using this,
%we will prove that $g$ has no symbolic extensions by means of Proposition \ref{prop4.4}.

%First, 
We define $\rho_0=\max\{\chi(\tilde{p},g); \; \tilde{p}\in Per_h(g) \text{ and related to } p_k(g)\}$, and  $$\xi_1(g)=\Big\{\mu_{\tilde{p}} : \tilde{p}\in Per_h(g), \text{ related to } p_k(g) \,\text{ and }\,\ \chi(\tilde{p},g)>{\displaystyle\frac{\rho_0}{2}}\Big\}$$
which is a non empty subset in $\SM(f)$.  Then, we consider the compact subset $\xi(g)=\overline{\xi_1(g)}$  in $\SM(g)$.
%Take $\frac{\chi(f)}{2}=\rho_0$ and

%To finished the proof   verify the conditions of Propositions \ref{propext} for the set $\xi_1$.  %Is sufficient prove the condition of the Prop \ref{propext} for $\mu=\mu_p \in \xi_1(f)$, since $\mu \in \xi_1(f)$, then $\xi(p)>\rho_0$.

%In oder to verify the hypothesis of Proposition \ref{prop4.4} for measures in $\xi$, and then conclude the proposition, it's enough to verify  only for
%measures in $\xi_1$.

Now, let $\mu_{\tilde{p}}\in \xi_1$ and $t$ be a positive integer. Since $g$ satisfies property $S_{n,p_k(g)}$ for every positive integer $n$,  there exist ergodic measures $\nu_m\rightarrow \mu_{\tilde{p}}$ such that $h_{\nu_m}(g)>\rho_0/2$, for every $m$. Moreover, %by property (b) of $S_{n,p_k(g)}$, for $m$ large enough, these measures are such that $h_t(\nu_m)=h_{\nu_m}(\alpha_t)=0$.
%Now, 
since these measures are supported on hyperbolic sets with the same index that $p_k(g)$,  by Sigmund \cite{S-I}, they are approximated by hyperbolic
periodic measures also supported in these hyperbolic sets, and by item (c) of property $S_{n,p_k(g)}$, they 
belong to $\xi_1(g)$. Hence, $\nu_m\in \xi(g)$ for every $m$, and then
$$
\limsup_{\nu_m\rightarrow \mu_{\tilde{p}}, \nu_m\in\xi(g)}h_{\nu_m}(g)>\frac{\rho_0}{2}.
$$

Therefore, since $p$ is arbitrary and $\xi_1(g)$ has dense periodic measures, there is no continuity point for the entropy function $h$. Thus, by Proposition \ref{propburguet}, this implies that $f$ has no symbolic extensions.

\begin{flushright}
 $\Box$
\end{flushright}

\section{Non existence of symbolic extensions versus good decomposition}\label{mainproof}

%\begin{proposition}\label{local}
% Let $f \in \mundo$, $p$ be a $f$-hyperbolic periodic point and $\SU(f)\subset\mundo$ be a small neighborhood of $f$. There is a residual subset
%$\mathfrak R \subset \SU(f)$ such that every $g \in \SU(f)$ satisfies only one of the following statements:
%\begin{itemize}
% \item[(i)] $H(p_g)$ has a dominated splitting;
%\item[(ii)] $g|_{H(p_g)}$ admits no symbolic extension.
%\end{itemize}
%Here, $H(p_g)$ denotes the continuation of $H(p,f)$.
%\end{proposition}

%%%%%%%%%%%%%%%%%%%%%%%%%% COMENTÁRIOS

%\begin{theorem}
%There is a residual subset $\mathfrak R \subset \mundo$ such that for all $f \in \mathfrak R$ and all $f$-hyperbolic periodic point $p$,
%only one of the following statements is true:
%\begin{itemize}
%\item[(i)] $H(p,f)$ has a good decomposition;
%\item[(ii)] $f|_{H(p_g)}$ has no symbolic extensions.
%\end{itemize}
%\end{theorem}

In this section we use Proposition \ref{local} and the generic machinery to prove Theorem \ref{maintheo}.
%The proof follows from known technics of generic theory. %the ideas of \cite{B-D}.

%Given $f\in\mundo$, recall that $Per^n(f)$ %=\{p\in M\;:\;f^j(p)=p,$ for some $j\leq n\}$
%is the set of all periodic points $p$ of $f$ such that the
%period of $p$ is smaller than or equal to $n$.

\vspace{0,2cm}
\textbf{ Proof of Theorem \ref{maintheo}:}

\vspace{0,1cm}

Since $\mundo$ is separable, there is a countable and dense subset $\mathcal{A}\subset \mundo$. Moreover, we can assume that $\mathcal{A}\subset \mathcal{R}_3$, the residual subset of $\mundo$ in the hypothesis of Proposition \ref{local}.

Now, for any $f\in \mathcal A$ and a small enough neighborhood $\mathcal{U}(f)$ of $f$, we consider the residual subset $\tilde{\mathcal{R}}_f$ in $\mathcal{U}(f)$ given by Proposition \ref{local}. Thus, we define
$$
\mathcal{R}_f=\tilde{\mathcal{R}}_f\cup (\mathcal{U}(f))^c,
$$
which is a residual subset in $\mundo$, indeed.  Also, since $\mathcal{A}$ is a dense subset %the following subset
$$
\mathcal{U}=\bigcup_{f\in \mathcal{A}} \mathcal{U}(f),
$$
is an open and dense subset of $\mundo$.

Finally, we define the following residual subset
$$
\mathcal{R}=\bigcap_{f\in \mathcal{A}} \mathcal{R}_f\cap \mathcal{U}.
$$

Now, let $g\in \mathcal{R}$  and $H(p,g)$ be a homoclinic class of $g$. Since $g\in \mathcal{U}$, there exists $f\in \mathcal{A}$ such that $g\in \mathcal{U}(f)$, and then provided $g$ also belongs to $\mathcal{R}_f$, $g$ should belongs to $\tilde{\mathcal{R}}_f$.  Therefore, by Proposition \ref{local} we have that either $H(p,g)$ has a good decomposition, or $f\,|\,H(p,g)$ has no symbolic extensions.  This completes the proof.

\section{The isolated case}\label{isolated}

It is enough to prove that for $f\in \SR$ of Theorem \ref{maintheo}, if $H(p,f)$ has a good decomposition and it is isolated then it is partially hyperbolic. Let $U$ be a neighborhood of $H(p,f)$ such that $H(p,f)=\bigcap_{n\in \Z}f^n(U)$. Also, let $E\oplus E_1\oplus\dots\oplus E_l\oplus F$ be the good decomposition. We will prove that $E$ is contracting, a similar argument will prove that $F$ is expanding.

We recall that $dim(E)$ is the smallest index of a periodic point in the class. By \cite{A-B-C-D-W}, there is another residual subset where we know that there exists a neighborhood $\SU$ of $f$ such that $dim(E)$ is still the smallest index of a periodic point in the class $H(p_g,g)$, for any $g\in \SU$, where $p_g$ is the analytic continuation of $p$. The intersection of this two residual subsets is the one claimed in the statement of the theorem.

Now, if $E$ does not contract, using the Ergodic Closing Lemma, as Ma\~n\'e did in \cite{M}, then it is possible to find $g\in \SU$ with a periodic orbit $O(q)\subset U$ with index smaller than $dim(E)$.

However, we can consider the result of Abdenur proved in \cite{Ab} to relative homoclinic class. Here, a relative homoclinic class of $p$ for $U$ is the subset of $H(p,f)$ of points that have the whole orbit inside $U$, which we denote by $H_U(p,f)$. Therefore, %there exists a residual subset of $\SU$ such that
since $H_U(p,f)=H(p,f)$ isolated homoclinic class, for any $h\in \mathcal{R}$ close enough to $f$,  $H_U(p(h),h)=\bigcap_{n\in \Z}h^n(U)$ and then $q\in H_U(p(h),h)\subset H(p(h),h)$, which is a contradiction.   %Obviously, we can assume that $g$ belongs to this residual subset. Hence $q\in H(p_g,g)$, which is a contradiction.

%Now we prove Corollary ?????

%Obviously, if the class has a principal symbolic extension then it has a symbolic extension. Now, if the class has a symbolic extension, then Theorem ????, the class is partially hyperbolic, with a good decomposition. Now, using [DIAZ....], the class has a principal symbolic extension.

\section{Proof of Proposition \ref{prophexp}}\label{hexpansivity}

First, note that one inclusion is a directly consequence of the result of Liao, Viana and Yang \cite{L-V-Y}. More precisely, they have proved that far away from homoclinic tangency every diffeomorphism is $h-$expansive.  For the other inclusion we will use Lemma \ref{mainlemma}.

Let $f\in HT$, that is, $f$ exhibits a homoclinic tangency, say $q$, for a hyperbolic periodic point $p$. Given $\epsilon > 0$ small, let us consider a small neighborhood $\mathcal{U}$  of $f$, with
$diam(\mathcal{U})<\epsilon$. Then by Lemma \ref{mainlemma}  there is a perturbation $f_1\in \mathcal{U}$ of $f$ such that $f_1$ has a periodic hyperbolic basic set $\La_1$ satisfying $$h(f_1|\La_1) >\chi(p(f_1),f_1)- \frac{1}{n_0+1},$$
for a big positive integer $n_0$ fixed, and moreover $f_1$ still exhibits a homoclinic tangency for $p(f_1)$.
 As we can see in the proof of the Lemma \ref{mainlemma}, $\Lambda_1$ can be found such that the base set $\overline{\Lambda}_1$, i.e., $\Lambda_1=\cup f_1^j(\overline{\Lambda}_1)$, is contained in a ball of radius so small, in particular, we can assume it is in a ball with radius $\frac{1}{n_0+1}$.

In the sequence, we consider a small neighborhood $\mathcal{U}_1$ of $f_1$, such that for all diffeomorphisms in $\mathcal{U}_1$  there is a continuation for $\Lambda_1$, and moreover $diam( \mathcal{U}_1)<\frac{\eps}{n_0+1}$. Now, using again Lemma \ref{mainlemma} we can find a diffeomorphism $f_2\in \mathcal{U}_1$, and a periodic hyperbolic basic set $\La_2$, with base set contained in a ball with radius $\frac{1}{n_0+2}$,  such that $$h(f_2|\La_2) >\chi(p(f_2),f_2)- \frac{1}{n_0+2},$$ and $f_2$ still exhibits a homoclinic tangency for $p(f_2)$.

Following this process inductively we can find a sequence of diffeomorphism $f_n\in \mathcal{U}_{n-1}$, with $diam( \mathcal{U}_n)<\frac{\eps}{n_0+n}$, $\mathcal{U}\supset \mathcal{U}_1\supset \ldots \supset \mathcal{U}_n\supset\ldots$, and moreover, by construction, $f_n$ is such that there exists periodic hyperbolic sets $\Lambda_1, \ldots, \Lambda_n$ with $diam( \Lambda_i)<\frac{1}{n_0+i}$, for every $1\leq i \leq n$, and
$$
h(f_n|\La_i) >\chi(p(f_n),f_n)- \frac{1}{n_0+i}.
$$
Since this sequence of diffeomorphism is a Cauchy sequence, it converges to a diffeomorphism $g$, that is $\epsilon$-close to $f$. Now, by choice of the open sets $\mathcal{U}_n$, $g$ has periodic hyperbolic basic sets with diameter so small as we want with topological entropy away from zero, since $\chi(p(f),f)$ varies continuously with the diffeomorphism $f$. Therefore, $g$ can not be asymptotically h-expansive. And then, we have proved that $\overline{HT}\subset \overline{NAHE}$.

\end{document}